\documentclass[a4paper, 12pt]{article}
\usepackage{cmap}					
\usepackage[T2A]{fontenc}			
\usepackage[utf8]{inputenc}			
\usepackage[english]{babel}
\usepackage{amsthm,amsmath,amsfonts,amssymb}
\usepackage{color}
\usepackage{graphicx,calc}
\usepackage{wrapfig}
\usepackage[font=small]{floatrow}
\usepackage{tikz}
\usepackage{floatflt}

\usepackage[hyper]{amsbib}

\usepackage[left=2.25cm,right=2.25cm,top=2.5cm,bottom=2.5cm,bindingoffset=0cm,nohead,nofoot,footskip=7mm]{geometry}

\newtheorem{theorem}{Theorem}           

\theoremstyle{definition}

\newtheorem{remark}{Remark}
\newtheorem{conj}{Conjecture}

\begin{document}

\title{Steklov-type 1D inequalities (a survey)}

\author{Alexander I. Nazarov and Alexandra P. Shcheglova}

%\address{Alexander I. Nazarov,\\ 
%St. Petersburg Dept. of Steklov Institute\\ 
%Fontanka 27, St. Petersburg, 191023, Russia\\ 
%and \\
%St. Petersburg State University,\\ 
%Universitetskii pr. 28, St. Petersburg, 198504, Russia \\
%\email{al.il.nazarov@gmail.com}}

%\address{Alexandra P. Shcheglova,\\ 
%St. Petersburg Electrotechnical University ``LETI''\\ 
%ul. Professora Popova 5, St. Petersburg, 197022, Russia\\
%and \\ 
%St. Petersburg State University,\\
%Universitetskii pr. 28, St. Petersburg, 198504, Russia \\
%\email{apshcheglova@etu.ru}}

\date{today}

%\subjclass{26D10}
        % AMS-2010 subj class. The list can be found on http://www.ams.org/mathscinet/msc/msc2010.html

%\thanks{We are grateful to Prof. N.G. Kuznetsov for valuable comments and suggestions, and to Prof. I.A. Sheipak who provided us with some important references}      

\maketitle

\begin{abstract}
We give a survey of classical and recent results on symmetry\,/\,asymmetry of extremal functions and sharp constants in $1$-dimensional functional inequalities.\medskip

{\bf Keywords}: one-dimensional functional inequalities;  symmetry; symmetry breaking; sharp constants.\medskip

{\bf MSC}: 26D10
\end{abstract}

\section{Introduction}

The first version of this survey was an extended variant of the talk given at the International Conference ``Qualitative Theory of Differential Equations'' in December 2020. In comparison with it, the present version is supplemented with both historical references and new results. 
The text has partial intersection with the survey \cite{KuzNaz} where an extensive bibliography on multidimensional functional inequalities is given.
\medskip

Sharp constants in one-dimensional functional inequalities are important in various fields of mathematics such as the theory of functions (see \cite{Maz}, \cite{Tikh}), mathematical physics (see, e.g., \cite{EK}, \cite{DELL}), mathematical statistics (see \cite[\S 6.2]{Nik}, \cite{LNN}) etc.

\medskip

The problem of sharp constants is closely related with the problem of symmetry and symmetry breaking of corresponding extremal functions. The first author was introduced to this topic almost 30 years ago by Yakov Yu. Nikitin and Vladimir A. Kondratiev. 

\medskip

The pioneer results of this type were sharp constants in the inequalities
\begin{align}
\label{eq:Steklov1}
&\int\limits_0^\ell u^2 (x) \, dx \le \Big( \frac{\ell}{\pi} \Big)^2 \int\limits_0^\ell [u'(x)]^2 \, dx, &&  \int\limits_0^\ell u(x)\,dx=0;\\ 
\label{eq:Steklov2}
&\int\limits_0^\ell u^2 (x) \, dx \le \Big( \frac{\ell}{\pi} \Big)^2 \int\limits_0^\ell [u'(x)]^2 \, dx, && u(0)=u(\ell)=0;\\ 
\label{eq:Steklov3}
&\int\limits_0^\ell u^2 (x) \, dx \le \Big( \frac{\ell}{2\pi} \Big)^2 \int\limits_0^\ell [u'(x)]^2 \, dx, && u(0)=u(\ell), \quad \int\limits_0^\ell u(x)\,dx=0.
\end{align}
These constants were found by Steklov \cite{St1} and \cite{St2}, respectively (see also \cite{St3}), and Almansi \cite{Al} (see also \cite{Kr}).\footnote{In the literature, inequalities (\ref{eq:Steklov1})--(\ref{eq:Steklov3}) are often called Poincar\'e inequalities or Wirtinger inequalities.} See \cite[Ch. II]{MPF}, \cite{KuzNaz} for a comprehensive history of inequalities (\ref{eq:Steklov1})--(\ref{eq:Steklov3}).

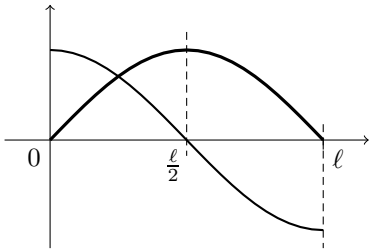
\begin{figure}[ht]
\begin{tikzpicture}[scale=2.0] 
\draw[->] (-0.5,0) -- (3.5,0);
\draw[->] (0,-1.2) -- (0,1.5);
\draw (3, -3pt) -- (3, 0) node[anchor=north west] {$\ell$};
\draw[densely dashed,thin] (1.5,1.2) -- (1.5,-5pt);
\draw (1.55,0pt) node[anchor=north east] {$\frac{\ell}{2}$} ;
\draw[densely dashed,thin] (3,5pt) -- (3,-1.2);
\draw[domain=0:3,color=black,very thick] plot (\x,{sin(60*\x)});
\draw[domain=0:3,color=black,thick] plot (\x,{cos(60*\x)});
\draw (0,0) node[anchor=north east] {$0$};
\end{tikzpicture} 
\vspace{4pt}
{\caption{The graphs of extremal functions for the inequalities (\ref{eq:Steklov1}) (thin line) and (\ref{eq:Steklov2}) (bold line).}\label{Fig1}}
\end{figure}

Notice that the sharp constants in (\ref{eq:Steklov1}) and (\ref{eq:Steklov2}) are attained by functions $\cos(\frac {\pi}{\ell}\,x)$ and $\sin(\frac {\pi}{\ell}\,x)$, respectively (up to a multiplicative constant). So, the extremal functions are symmetric (respectively, odd and even) about the middle of the interval; see Fig.~\ref{Fig1}. The extremal functions in (\ref{eq:Steklov3}) are given by any linear combination of 
$\cos(\frac {2\pi}{\ell}\,x)$ and $\sin(\frac {2\pi}{\ell}\,x)$.

\section{A simplest extension of the inequality (\ref{eq:Steklov2})}

We begin with the inequality
\begin{equation}
    \label{eq:Friedrichs}
    \|u\|_{L^q(0,\ell)}\le \lambda_1\ell^{\frac 1q+\frac 1{p'}}\|u'\|_{L^p(0,\ell)}, \qquad u(0)=u(\ell)=0.
\end{equation}
(here and below $1\le p,q\le\infty$, and $p'$ stands for the H\"older conjugate exponent to $p$). By dilation, it is easy to see that $\lambda_1$ depends only on $p$ and $q$. So, it is sufficient to consider $\ell=1$.
\medskip

Inequality (\ref{eq:Friedrichs}) is equivalent to several ones. We list three examples:
\begin{equation}
\label{eq:Friedrichs1}
\aligned
&  \|u\|_{L^q(0,1)}\le 2\lambda_1\|u'\|_{L^p(0,1)}, && u(0)=0;\\
&  \|u\|_{L^q(0,1)}\le \lambda_1\|u'\|_{L^p(0,1)}, && u(0)+u(1)=0;\\
&  \|u\|_{L^q(0,1)}\le \frac 12\lambda_1\|u'\|_{L^p(0,1)}, && u(0)=u(1),\quad \min u+\max u=0.
\endaligned
\end{equation}

\begin{figure}[ht]
\begin{tikzpicture}[scale=2.2] 
\draw[->] (-0.5,0) -- (3.5,0);
\draw[->] (0,-0.2) -- (0,2);
\draw (3, -3pt) -- (3, 0) node[anchor=north west] {$1$};
\draw[densely dashed,thin] (1.5,2) -- (1.5,-5pt);
\draw (1.55,0pt) node[anchor=north east] {$\frac{1}{2}$} ;
%\draw[densely dashed,thin] (3,5pt) -- (3,-1.2);
\draw[domain=0:3,color=black, thick] plot (\x,{0.8*\x*(3-\x)});
\draw[very thick] (0,0) -- (1.5,1.5) -- (3,0);
\draw (0,0) node[anchor=north east] {$0$};
\end{tikzpicture} 
\vspace{4pt}
{\caption{The graphs of extremal functions in (\ref{eq:Friedrichs}) for $p=2$: $q=1$ (thin line) and $q=\infty$ (bold line).}\label{Fig2}}
\end{figure}
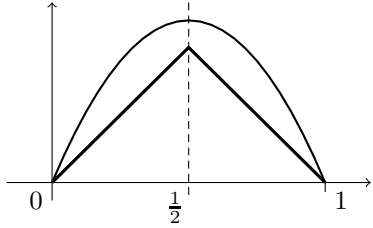

The following statement holds:

\begin{theorem}\label{T1}
The sharp constant in (\ref{eq:Friedrichs}) is given by
\begin{equation*}
\label{eq:Schmidt}
\lambda_1 (p,q) = \frac{\mathfrak F \big( \frac 1{p'} + \frac 1{q} \big)}{2 \, \mathfrak F
\big( \frac 1{p'} \big) \mathfrak F \big( \frac 1{q} \big)} \, , 
\end{equation*}
where $\mathfrak F (s) = \frac{\Gamma (s+1)}{s^s}$. The corresponding extremal function $U_{p,q}$ can be expressed in quadratures, does not change sign and is even with respect to $x=\frac 12$, see Fig.~\ref{Fig2}.
\end{theorem}

\begin{remark}
For $p>1$, the natural domain for $u$ in (\ref{eq:Friedrichs}) is the Sobolev space 
$\stackrel{\circ\ }{W^1_p}\!\!(0,1)$. For $p=1$ the sharp constant in (\ref{eq:Friedrichs}) is not achieved in the Sobolev space. However, in this case one can consider $u\in BV(0,1)$ and understand the right-hand side of (\ref{eq:Friedrichs}) in the sense of measures. In this case, the statement of Theorem \ref{T1} is true (except for the case $p=1$, $q=\infty$, where there are symmetric and asymmetric extremals).
\end{remark}

%\medskip

The history of Theorem \ref{T1} is given in the Table \ref{tab1}.\footnote{In fact, Hardy--Littlewood--P\'olya and Levin dealt with the first inequality in (\ref{eq:Friedrichs1}) whereas Schmidt considered the third inequality in (\ref{eq:Friedrichs1}).}
\medskip

\begin{table}[ht]    
\centering
\begin{tabular} {|c|c|l|}
\hline Year & Authors & Parameters \\
\hline 1901 & Steklov \cite{St2} & $p=q=2$ \\
\hline 1934 & Hardy et al \cite[Sec. 7.6]{HLP} & $p=q=2k$, $k\in\mathbb N$ \\
\hline 1938 & Levin \cite{Le} & $p=q$ \\
\hline 1940  & Schmidt \cite{Schm} & $\forall\, p,q$ \\
\hline\end{tabular}
    \caption{The history of inequality (\ref{eq:Friedrichs})}
    \label{tab1}
\end{table}

\begin{remark}
Notice that the result of Theorem \ref{T1} and its particular cases were later rediscovered several times even in the 21st century; see, for instance, \cite[p. 50--51]{Bee}, \cite[p. 220]{Str}, \cite[p. 377]{Bo}, \cite[p. 357]{Tal}, \cite[Th.~5.1]{DraMan}, \cite{BS}, \cite{EL}.
\end{remark}

\section{An extension of the inequality (\ref{eq:Steklov1})}

Now we consider the following inequality: 
\begin{equation}
    \label{eq:Poincare}
    \|u\|_{L^q(0,\ell)}\le \lambda_2\ell^{\frac 1q+\frac 1{p'}}\|u'\|_{L^p(0,\ell)}, \qquad \int\limits_0^\ell |u(x)|^{r-2}u(x)\,dx=0
\end{equation}
(here $1\le p,q,r\le\infty$; for $r=\infty$ the last relation is understood in the limit sense). As in (\ref{eq:Friedrichs}), $\lambda_2=\lambda_2(p,q,r)$ does not depend on $\ell$, and we put $\ell=1$. 
\medskip

Under additional restriction that $u$ is $1$-periodic, inequality (\ref{eq:Poincare}) holds with the sharp constant $\frac 12 \lambda_2$. In the case $r=2$, (\ref{eq:Poincare}) is equivalent to several other inequalities. We again list three examples:
{\allowdisplaybreaks
\begin{align}
\label{eq:Poincare1}
& \|u'\|_{L^q(0,1)}\le \frac 12\lambda_2\|u''\|_{L^p(0,1)},
&& u(0)=u(1)=u'(0)=u'(1)=0;\\
\nonumber
& \|u^{(k)}\|_{L^q(0,1)}\le \frac 12\lambda_2\|u^{(k+1)}\|_{L^p(0,1)},
&& u\ \mbox {is $1$-periodic},\quad  k\in\mathbb N;\\
\nonumber
& \|u-\overline u\|_{L^q(0,1)}\le \frac 12\lambda_2\|u'\|_{L^p(0,1)},
&& u(0)=u(1)=0
\end{align}}
($\overline u$ stands for the mean value of $u$).
\bigskip

Notice that if the extremal function in (\ref{eq:Poincare}) is odd w.r.t. $x=\frac 12$ then the integral restriction is fulfilled for any $r$. Therefore, in this case, $\lambda_2$ does not depend on $r$. However, the general picture is more complicated.
\medskip

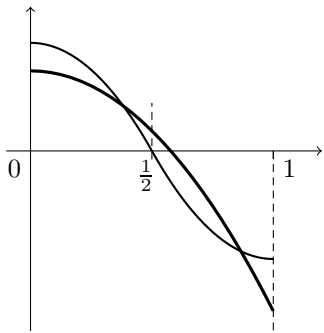
\begin{figure}[ht]
\begin{tikzpicture}[scale=4.3] 
\draw[->] (-0.1,0) -- (1.2,0);
\draw[->] (0,-0.75) -- (0,0.6);
\draw (1, -1pt) -- (1, 0) node[anchor=north west] {$1$};
\draw[densely dashed,thin] (1,0) -- (1,-0.75);
\draw[densely dashed,thin] (0.5,-1pt) -- (0.5,0.2);
\draw (0.55,0pt) node[anchor=north east] {$\frac{1}{2}$} ;
%\draw[densely dashed,thin] (3,5pt) -- (3,-1.2);
\draw[domain=0:1,color=black,very thick] plot (\x,{1/3-\x*\x});
\draw[domain=0:0.5,color=black,thick] plot (\x,{1.8*(1/4-\x*\x)});
\draw[domain=0.5:1,color=black,thick] plot (\x,{1.8*(\x*\x-2*\x+3/4)});
\draw (0,0) node[anchor=north east] {$0$};
\end{tikzpicture} 
\vspace{4pt}
{\caption{The graphs of extremal functions in (\ref{eq:Poincare}) for $p=2$, $r=2$: $q=1$ (thin line) and $q=\infty$ (bold line). }\label{Fig3}}
\end{figure}

\begin{theorem}\label{T2}
If $q\le (2r-1)p$ then the following equality holds:
$$
\lambda_2(p,q,r)=\lambda_1(p,q),
$$
see (\ref{eq:Schmidt}). The corresponding extremal function $V_{p,q}$ is (up to a multiplicative constant) given by formula
$$
V_{p,q}(x) = \begin{cases} 
\phantom{-\,}
U_{p,q} \big( x+\frac{1}{2} \big) & {\rm if}\ x\le \frac{1}{2} , \\
-\, U_{p,q} \big( x-\frac{1}{2} \big) & {\rm if}\ x \geq
\frac{1}{2} , 
\end{cases}
$$
where $U_{p,q}$ is introduced in Theorem \ref{T1}. In particular, $V_{p,q}$ is odd w.r.t. $x=\frac 12$.\medskip

In contrast, if $q> (2r-1)p$ then $\lambda_2(p,q,r)>\lambda_1(p,q)$, and the extremal function $V_{p,q}$ has no symmetry; see Fig.~\ref{Fig3}.
\end{theorem} 

\begin{remark}
Similarly to Theorem \ref{T1}, for $p=1$ the statement of Theorem \ref{T2} is true if one understands the right-hand side of (\ref{eq:Poincare}) in the sense of measures.

If $r=1$ then the statement is true for $q\le p$ whereas for $q>p$ the sharp constant in (\ref{eq:Poincare}) is not achieved.\footnote{However, \ for $q>p$ any normalized minimizing sequence converges to a non-symmetric function. See \cite[\S 2]{GeNaz};  cf. also \cite{BT}.}

For $r=\infty$ the last relation in (\ref{eq:Poincare}) should be understood in the limit sense, that is $\min u+\max u=0$. The statement of Theorem \ref{T2} is true.
\end{remark}
\medskip

The history of Theorem \ref{T2} is given in Table \ref{tab2}.\footnote{In fact, Bohr, Dacorogna--Gangbo--Subia, Belloni--Kawohl and Croce--Dacorogna dealt with (\ref{eq:Poincare}) for periodic functions, while Egorov and Buslaev--Kondratiev--Nazarov considered inequality (\ref{eq:Poincare1}). }

\begin{remark}
In \cite[Sect{.}~1.1.19]{Maz}, the equality $\lambda_2(p,p,2)=\lambda_1(p,p)$ is attributed to A.~Stanoyevitch \cite{Sta}. However, the proof in \cite{Sta} turned out to be incorrect.

A part of Theorem \ref{T2} proved in \cite{DGS} was also rediscovered later. See \cite{EHS} for the case $p=q=r$ and \cite{BS1} for $r=2$, $q\le 2p$. On the other hand, it is erroneously stated in \cite[Th.~5.2]{DraMan} that 
$\lambda_2(p,q,2)=\lambda_1(p,q)$ for any $q$.
\end{remark}

\begin{table}[ht]    
\centering
\begin{tabular} {|c|c|l|l|}
\hline Year & Authors & Symmetry & Asymmetry \\
\hline 1896 & Steklov \cite{St1} & $r=2$, $q=p=2$ &  \\
\hline 1935 & Bohr \cite{Bohr} & $r=2$, $q=p=\infty$ &  \\
\hline 1992 & Dacorogna et al \cite{DGS} & $r=q$ &  \\
\hline 1992 & \cite{DGS} & $r=2$, $q\le 2p$ & $r=2$, $q>>1$ \\
\hline 1997 & Egorov \cite{E} & & $r=2$, $q>4p-1$ \\
\hline 1998 & Buslaev et al \cite{BKN} & $r=2$, $q\le2p+\varepsilon$ & $r=2$, $q>3p$ \\
\hline 1999 & Belloni, Kawohl \cite{BeKa}\footnotemark & $r=2$, $q\le 2p+1$ &  \\
\hline 2002 & Nazarov \cite{Naz_2002} & $r=2$, $q\le 3p$ &  \\
\hline 2002 & Abessolo \cite{Ab} & $q\le rp+\varepsilon$ & $q>r^2p-(r-1)^2$ \\
\hline 2003 & Croce, Dacorogna \cite{CrDa} & $q\le rp+r-1$ & $q>(2r-1)p$ \\
\hline 2011 & Gerasimov, Nazarov \cite{GeNaz}\footnotemark[6] & $q\le (2r-1)p$ & \\
\hline 2018 & Ghisi et al \cite{GGR}\footnotemark[6] & $q\le (2r-1)p$ &  \\
\hline\end{tabular}
    \caption{The history of inequality (\ref{eq:Poincare})}
    \label{tab2}
\end{table}
\footnotetext{A technical gap in the proof was fixed in \cite{Ka}.}
\addtocounter{footnote}{1}
\footnotetext{In \cite{GeNaz}, a computer-assisted proof was given, while Ghisi--Gobbino--Rovellini succeeded in a pure analytical proof. }

%\medskip
Up to our knowledge, a unique explicit expression for $\lambda_2$ in the asymmetry region of parameters is as follows.

\begin{theorem}[{\cite[\S 4]{BKN}}; see also \cite{Stech}]\!\!\footnote{Stechkin considered (\ref{eq:Poincare}) and some higher-order inequalities for periodic functions.} \ 
$\lambda_2(p,\infty, 2)=(p'+1)^{-\frac 1{p'}}$.
\end{theorem}

\section{A higher-order extension of inequality (\ref{eq:Steklov2})}

Further, we consider the following inequality: 
\begin{equation}
\label{eq:high}
\|u^{(k)}\|_{L^q(0,\ell)}\le \lambda_3\ell^{n-k+\frac 1q-\frac 1p}\|u^{(n)}\|_{L^p(0,\ell)}, \qquad u\in\stackrel{\circ\ }{W^n_p}\!\!(0,\ell)
\end{equation}
(here $n,k\in\mathbb Z_+$, $n>k$). 

\begin{remark}
As earlier, $\lambda_3=\lambda_3(n,k,p,q)$ does not depend on $\ell$, and we can put $\ell=1$. For $p=1$ the right-hand side of (\ref{eq:high}) is understood in the sense of measures. 
Evidently, $\lambda_3(1,0,p,q)=\lambda_1(p,q)$. Moreover, by (\ref{eq:Poincare1}) we have $\lambda_3(2,1,p,q)=\frac 12 \lambda_2(p,q,2)$.
\end{remark}

\begin{remark}
Recall that for $u\in\stackrel{\circ\ }{W^n_p}\!\!(0,\ell)$, {\bf all} derivatives of $u$ up to $(n-1)$-th order vanish at the endpoints of the interval. In this connection we mention the paper \cite{BoLa}, where the inequality
\begin{equation*}
\|u\|_{L^q(0,\ell)}\le \widehat\lambda_3\ell^{2+\frac 1q-\frac 1p}\|u''\|_{L^p(0,\ell)}, \qquad u(0)=u(\ell)=0,
\end{equation*}
was investigated. In particular, it was proved that the extremal function in this inequality is even w.r.t. the middle of the interval. Also in the case $q=p'$ the sharp constant in this inequality was (somewhat unexpectedly) calculated explicitly for general $p$.\footnote{Evidently, $\widehat\lambda_3(p,q)>\lambda_3(2,0,p,q)$ for all $p,q$.}
It turns out that $\widehat\lambda_3(p,p')= \lambda_1^2(2,p')$, and the extremal function is just $U_{2,p'}$, see Theorem~\ref{T1}.
\end{remark}

%\medskip

The results listed earlier, as well as some calculations, lead to the following conjecture; see \cite{MukNaz}.\medskip

\begin{conj}
 If $k$ is even, then the extremal function in the problem (\ref{eq:high}) is even w.r.t. $x=\frac 12$ for all admissible $n, p, q$ (except for the case\footnote{In this case extremal function can be as symmetric as asymmetric; cf. \cite[Th.~3]{GarSh_24}, where it was shown that $\lambda_3(n,n-1,1,\infty)=\frac 12$.} $p=1$, $q=\infty$, $n=k+1$). If $k$ is odd, then for all admissible $n$ and $p$ there exists $\widehat q(n,k,p)>p$ such that the extremal is even w.r.t. $x=\frac 12$ for $q\le \widehat q$ and is non-symmetric for $q>\widehat q$.
\end{conj}

The known values of $\lambda_3$, besides $n=1$, $k=0$, and  $n=2$, $k=1$, concern the cases $p=2$ ($q=1,2$) or $q=\infty$.\footnote{In \cite{Kal_14}, two-sided estimates of $\lambda_3(n,0,2,q)$ were obtained for general $q$.}
\medskip

The following statement was proved in \cite{J3}\footnote{The case $n=2$ was considered in \cite{J1}. The announcement without proof was given in \cite{J2}. Later, the result of \cite{J3} was rediscovered in \cite{NazPet}.} for $k=n-1$ and in \cite{Pet}\footnote{See also \cite{Cimm} and \cite{J4} (without proof), and \cite{Sl} for $k=n-2$.} in the general case.

\begin{figure}[ht]
\begin{tikzpicture}[scale=4.5] 
\draw[->] (-0.1,0) -- (1.1,0);
\draw[->] (0,-0.1) -- (0,1);
\draw (0.56,0pt) node[anchor=north east] {$\frac{1}{2}$} ;
\draw[densely dashed,thin] (0.5,-1pt) -- (0.5,0.95);
\draw (0,0) node[anchor=north east] {$0$};
\draw (1,-1pt) -- (1,0) node[anchor=north west] {$1$};
\draw[yscale=1.7,very thick] plot[smooth,color=black] file {fig4_3122.table};
\draw[yscale=1.7, thick] plot[smooth,color=black] file {fig4_4222.table};
\end{tikzpicture} 
\vspace{4pt}
{\caption{The graphs of extremal functions in (\ref{eq:high}) for $p=2$, $q=2$: $n=4$, $k=2$ (thin line) and $n=3$, $k=1$ (bold line).} \label{Fig4}}
\end{figure}
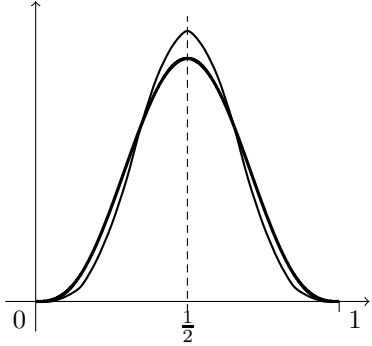

\begin{theorem} 
\label{T3}
$\omega=\lambda_3^{-\frac 1{n-k}}(n,k,2,2)$ is the least positive root of the function
$$
\Phi^{(n,k)}(\omega)=\det \big[\mathcal{D}^{(n,k)}(\omega)\big],
$$
where $\mathcal{D}^{(n,k)}(\omega)$ is the $(n-k)\times(n-k)$-matrix with entries
$$
\mathcal{D}^{(n,k)}_{jm}(\omega)=
(\omega z^m)^{\frac {2k+2j-1}2} J_{\frac {2k+2j-1}2}(\omega z^m), \quad  j,m=0,\dots, n-k-1
$$
(here $z=e^{\frac {i\pi}{n-k}}$, whereas
$J_\nu$ is the Bessel function of the first kind).
The corresponding extremal function is even w.r.t. $x=\frac 12$, see Fig.~\ref{Fig4}.
\end{theorem}

For $q=\infty$ and general $p$, the answers are known only for $n=2$, $k=0$ \cite{Osh, WKNY1} and $n=3$, $k=0$ \cite{WKNY1}.\footnote{Recently Garmanova--Sheipak \cite{GarSh_24} derived a general relation between 
$\lambda_3(n,k,p,\infty)$ and the best approximation of a special spline by polynomials in $L^{p'}(0,1)$. However, at the moment the explicit result is obtained only in some particular cases, see below.}

\begin{theorem} 
\label{T4}
The following equalities hold:
\begin{eqnarray*}
\lambda_3(2,0,p,\infty) &=& \frac 18\,(p'+1)^{-\frac 1{p'}}.\\
\lambda_3(3,0,p,\infty) &=& \frac 1{16} \cdot\min\limits_{\alpha\in(0,1)} \Big(\int\limits_0^1 x^{p'} |x-\alpha|^{p'}\,dx\Big)^{\frac 1{p'}}.
\end{eqnarray*}
The corresponding extremal function is even w.r.t. $x=\frac 12$.
\end{theorem}

It is convenient to introduce the function
$$
A_{n,k,p}(a)=\max\{|u^{(k)}(a)|\,:\, u\in\stackrel{\circ\ }{W^n_p}\!\!(0,1),\ \|u^{(n)}\|_{L^p(0,1)}\le 1\}, \quad a\in(0,1).
$$

\begin{remark}
\label{remark7}
It is evident that $\lambda_3(n,k,p,\infty)=\max\limits_{a\in(0,1)}A_{n,k,p}(a)$. Moreover, it is easy to see that the extremal function in (\ref{eq:high}) with $1<p<\infty$ and $q=\infty$ is even w.r.t. $x=\frac 12$ if and only if $\max\limits_{a\in(0,1)}A_{n,k,p}(a)=A_{n,k,p}(\tfrac 12)$.
\end{remark}

\begin{remark}
The function $A_{n,k,2}(a)$ was first introduced in \cite{Kal}. In particular, it was shown in this paper that $A_{n,k,2}^2(a)$ is a degree $2n-1$ polynomial of the variable $t=a-a^2\in(0,\frac 14)$.
\end{remark}

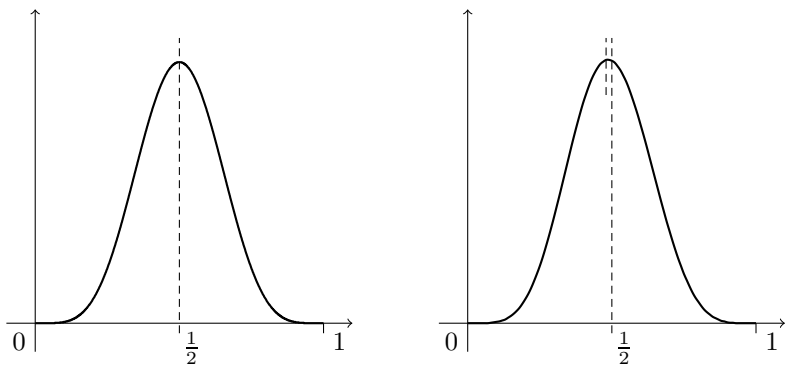
\begin{figure}[ht]
\begin{tikzpicture}[xscale=5.5, yscale=4.5]
\draw[->] (-0.1,0) -- (1.1,0);
\draw[->] (0,-0.1) -- (0,1.1);
\draw (0.48,0pt) node[anchor=north west] {$\frac{1}{2}$} ;
\draw[densely dashed,thin] (0.5,-1pt) -- (0.5,1);
\draw (0,0) node[anchor=north east] {$0$};
\draw (1,-1pt) -- (1,0) node[anchor=north west] {$1$};
\draw[thick] plot[yscale=1.5,smooth,color=black] file {fig5_52.table};

\begin{scope}[xshift=1.5cm]
\draw[->] (-0.1,0) -- (1.1,0);
\draw[->] (0,-0.1) -- (0,1.1);
\draw (0.48,0pt) node[anchor=north west] {$\frac{1}{2}$} ;
\draw[densely dashed,thin] (0.5,-1pt) -- (0.5,1);
\draw[densely dashed,thin] (0.48,0.8) -- (0.48,1);
%\draw[densely dashed,thin] (0.48,-1pt) -- (0.48,1);
\draw (0,0) node[anchor=north east] {$0$};
\draw (1,-1pt) -- (1,0) node[anchor=north west] {$1$};
\draw[thick] plot[yscale=1.6, smooth,color=black] file {fig5_53.table};
\end{scope}
\end{tikzpicture} 
\vspace{4pt}
{\caption{The graphs of extremal functions in (\ref{eq:high}) for $p=2$, $q=\infty$: $n=5$, $k=2$ (left) and  $n=5$,~$k=3$ (right).} \label{Fig5}}
\end{figure}

The following statement holds:

\begin{theorem} 
\label{T5}
Let $k$ be even. Then
$$
\lambda_3(n,k,2,\infty)=A_{n,k,2}(\tfrac 12)=\frac {(k-1)!!}{2^{2n-\frac {3k}2-1} (n-\frac k2-1)!\sqrt{2n-2k-1}}.
$$
The corresponding extremal function is even w.r.t. $x=\frac 12$.
\medskip

In contrast, if $k$ is odd, then $\lambda_3(n,k,2,\infty)>A_{n,k,2}(\tfrac 12)$, and the extremal function has no symmetry; see Fig.~\ref{Fig5}.
\end{theorem}

If $k$ is odd then explicit expressions for $\lambda_3(n,k,2,\infty)$ are known only for $k=1$ \cite{Kal}, $k=3,5$ \cite{ShGar1}.\footnote{In \cite{Gar_21}, two-sided estimates of $\lambda_3(n,k,2,\infty)$ were obtained for odd $k$.}

\begin{theorem} 
\label{T6}
The following equalities hold:
$$
\lambda_3(n,1,2,\infty)=A_{n,1,2}(a_{1,2})=\frac 1{(n-1)!}\Big(\frac  {n-1}{2(2n-1)}\Big)^{n-1}\sqrt{\frac {2n-1}{2n-3}},
$$
where $a_{1,2}=\frac 12\big(1\mp\frac 1{\sqrt{2n-1}}\big)$;

\begin{multline*}
\lambda_3(n,3,2,\infty)=A_{n,3,2}(a_{3,4})=\Big(\tfrac {(n-2)(2n-3)+\sqrt{(n-2)(2n-3)}}{2(2n-1)(2n-3)}\Big)^{n-\frac 72}\\
\times\tfrac {\sqrt{3(n-2)(2n-5)-(2n-7)\sqrt{3(n-2)(2n-3)}}}  {(n-2)!(2n-1)}\sqrt{\tfrac {n-2}{2(2n-7)}},
\end{multline*}
where 
$$
a_{3,4}=\tfrac 12\Big(1\mp\sqrt{1-2\,\tfrac {(n-2)(2n-3)+\sqrt{3(n-2)(2n-3)}}{(2n-1)(2n-3)}}\Big);
$$

$$
\lambda_3(n,5,2,\infty)=A_{n,5,2}(a_{5,6})\equiv A_{n,5,2}\big(\tfrac 12\big(1\mp\sqrt{1-4{\bf t}} \big)\big),
$$
where
$$
{\bf t}=\tfrac {n-3}{2(2n-1)}+\tfrac {\sqrt{5(n-3)}}{(2n-1)\sqrt{2n-3}}\cos\Big(\tfrac 13 \arccos\Big(-\tfrac {2n-11}{2n-5}\sqrt{\tfrac{2n-3}{5(n-3)}}\Big)\Big).
$$
The explicit expression for $\lambda_3(n,5,2,\infty)$ is quite complicated, and we omit it.
\end{theorem}

Figure~\ref{Fig6} shows examples of functions $A_{n,k,2}(a)$ for even and odd $k$.
%\footnote{The scales of the left and right pictures are different.}

\begin{figure}[ht]
\begin{tikzpicture}[scale=4.0] 
\draw[->] (-0.1,0) -- (1.1,0);
\draw[->] (0,-0.1) -- (0,0.9);
\draw (0.55,0pt) node[anchor=north east] {$\frac{1}{2}$} ;
\draw[densely dashed,thin] (0.5,-1pt) -- (0.5,0.9);
\draw (0,0) node[anchor=north east] {$0$};
\draw (1,-1pt) -- (1,0) node[anchor=north west] {$1$};
\draw[thick] plot[smooth,color=black] file {fig6_A_64.table};

\begin{scope}[xshift=1.5cm]
\draw[->] (-0.1,0) -- (1.1,0);
\draw[->] (0,-0.1) -- (0,0.9);
\draw (0.55,0pt) node[anchor=north east] {$\frac{1}{2}$} ;
\draw[densely dashed,thin] (0.5,-1pt) -- (0.5,0.6);
\draw[densely dashed,thin] (0.375,0.84) -- (0.375,-1pt) node[anchor=north] {$a_3$} ;
\draw[densely dashed,thin] (0.625,0.84) -- (0.625,-1pt) node[anchor=north] {$a_4$} ;
\draw (0,0) node[anchor=north east] {$0$};
\draw (1,-1pt) -- (1,0) node[anchor=north west] {$1$};
\draw[thick] plot[smooth,color=black] file {fig6_A_63.table};
\end{scope}
\end{tikzpicture} 
\vspace{4pt}
{\caption{The graphs of the functions $A_{n,k,2}(a)$ for $n=6$: $k=4$ (left) and $k=3$ (right). The scales of the pictures are different.} \label{Fig6}}
\end{figure}

For $p=q=\infty$, the sharp constant is known in the cases $k=n-1$ \cite{GarSh_24} and $k=0$ \cite{KazSh}.

\begin{theorem}
\label{T7}
\begin{enumerate}
    \item 
Let $k=n-1$ be even. Then    
$$
\lambda_3(n,n-1,\infty,\infty)=A_{n,n-1,\infty} \big(\tfrac 12\big) =\frac {\tan\big(\tfrac{\pi}{2(n+1)}\big)} 2.
$$
The corresponding extremal function is even w.r.t. $x=\frac 12$.

\item
Let $k=n-1$ be odd. Then    
$$
\aligned
\lambda_3(n,n-1,\infty,\infty)= &\, A_{n,n-1,\infty }(\mathfrak{a}_{\frac n2})
\\
= &\, A_{n,n-1,\infty}(\mathfrak{a}_{\frac n2+1}) = \frac{\tan\big(\tfrac{\pi}{2(n+1)}\big)\sin\big(\tfrac{\pi n}{2(n+1)}\big)} 2,
\endaligned
$$
where $\mathfrak{a}_j=\sin^2\tfrac{\pi j}{4(n+1)}$, $j=1,\dots, n$.

The corresponding extremal function is even w.r.t. $x=\frac 12$, see Fig. \ref{Fig7}.

\end{enumerate}

\end{theorem}

\begin{figure}[ht]
\begin{tikzpicture}[xscale=5,yscale=3.5] 
\draw[->] (-0.1,0) -- (1.1,0);
\draw[->] (0,-0.1) -- (0,1.15);
\draw (0.56,0pt) node[anchor=north east] {$\frac{1}{2}$} ;
\draw[densely dashed,thin] (0.5,-1pt) -- (0.5,1.1);
\draw (0,0) node[anchor=north east] {$0$};
\draw (1,-1pt) -- (1,0) node[anchor=north west] {$1$};
\draw[yscale=1.7,very thick] plot[smooth,color=black] file {fig8_43.table};
\draw[yscale=1.7, thick] plot[smooth,color=black] file {fig8_32.table};
\end{tikzpicture} 
\vspace{4pt}
{\caption{The graphs of extremal functions in (\ref{eq:high}) for $p=q=\infty$: $n=3$, $k=2$ (thin line) and $n=4$, $k=3$ (bold line).} \label{Fig7}}
\end{figure}

Figure~\ref{Fig8} shows examples of functions $A_{n,n-1,\infty}(a)$ for odd and even $n$. The algorithm for constructing these graphs was given in the paper~\cite{Dei}.

\begin{figure}[ht]
\begin{tikzpicture}[xscale=4.0, yscale=8.0] 
\draw[->] (-0.1,0) -- (1.1,0);
\draw[->] (0,-0.1) -- (0,0.3);
\draw[densely dashed,thin] (0.5,0.3) -- (0.5,0) node[anchor=north] {$\frac{1}{2}$};
\draw (0,0) node[anchor=north east] {$0$};
\draw (1,0) node[anchor=north west] {$1$};
\draw[thick] (0,0) to (0.146,0.146) to[out=0,in=-160] (0.5,0.207) to[out=-20,in=180] (0.854,0.146) to (1,0);
\draw[densely dotted,thin] plot[smooth,color=black!70] file {fig7_B3.table};

\begin{scope}[xshift=1.5cm, yscale=1.25]
\draw[->] (-0.1,0) -- (1.1,0);
\draw[->] (0,-0.1) -- (0,0.24);
\draw[densely dashed,thin] (0.5,0.24) -- (0.5,0) node[anchor=north] {$\frac{1}{2}$};
\draw[densely dashed,thin] (0.345,0.154) -- (0.345,-0.01) node[anchor=north] {$\mathfrak{a}_2$} ;
\draw[densely dashed,thin] (0.655,0.154) -- (0.655,-0.01) node[anchor=north] {$\mathfrak{a}_3$} ;
\draw (0,0) node[anchor=north east] {$0$};
\draw (1,0) node[anchor=north west] {$1$};
\draw[thick] (0,0) to (0.095,0.095) to[out=0,in=-160] (0.345,0.154) to[out=-10,in=-170]  (0.655,0.154) to[out=-20,in=180] (0.905,0.095) to (1,0);
\draw[densely dotted,thin] plot[smooth,color=black!70] file {fig7_B4.table};
\end{scope}
\end{tikzpicture} 
\vspace{4pt}
{\caption{The graphs of the functions $A_{n,n-1,\infty}(a)$ for $n=3$ (left) and $n=4$ (right). The graphs of the functions $B_n(a)=\tan\big(\tfrac{\pi}{2(n+1)}\big)\sqrt{a-a^2}$ are given by the dotted line. The scales of the pictures are different.
} 
\label{Fig8}}
\end{figure}
%Обозначения

\begin{remark}
Notice that for odd $k=n-1$ the extremal function is even, despite the fact that
$\max\limits_{a\in(0,1)}A_{n,n-1,\infty}(a)>A_{n,n-1,\infty}(\tfrac 12)$%, see Fig. \ref{Fig8}
. This distinguishes the case $p=\infty$ from the case $p<\infty$, cf. Remark \ref{remark7}.

The symmetry of the extremal function was not discussed in \cite{GarSh_24}. We prove it in the Appendix.
\end{remark}

\begin{theorem}
\label{T8}
The following equality holds:
$$
\lambda_3(n,0,\infty,\infty)=A_{n,0,\infty}(\tfrac 12)=\frac{n+1}{2^{2n-2}\pi n!}\int\limits_{0}^1  \frac{(1-x^2)^n}{1+(-1)^n x^{2(n+1)}}\, dx.
$$
This integral can be expressed in terms of hypergeometric functions, see \cite{KazSh}.
The corresponding extremal function is even w.r.t. $x=\frac 12$.
\end{theorem}
%The paper \cite{KazSh} also gives an expression for this integral in terms of the hypergeometric functions.
\medskip

We turn to the case $q=1$. The first result here was obtained quite recently \cite{HNOR}.

\begin{theorem}
\label{T9}
The following equality holds:
$$
\lambda_3(n,0,2,1)=\frac{n!}{(2n)!\sqrt{2n+1}}.
$$
The corresponding extremal function equals $x^n(1-x)^n$ (up to a multiplicative constant). Evidently, it is even w.r.t. $x=\frac 12$.
\end{theorem}

%\bigskip

The history of Theorems \ref{T3}--\ref{T9} is given in Table \ref{tab3}.

\begin{table}[ht]    
\centering
    \begin{tabular} {|c|c||c|c|c||l|l|}
\hline Year & Authors & $n$ & $k$ & $p$ & Symm. & Asymm. \\
\hline 1940 & Schmidt \cite{Schm} & $1$ & $0$ & $\forall$ & $\forall q$ &  \\
\hline 1998 & Buslaev et al \cite{BKN} & $2$ & $1$ & $\forall$ &  &  $q>3p$ \\
\hline 2002 & Nazarov \cite{Naz_2002} & $2$ & $1$ & $\forall$ & $q\le 3p$ &  \\
\hline 1931 & Janet \cite{J3}\footnotemark & $\forall$ & $n-1$ & $2$ & $q=2$ &  \\
\hline 2017 & Yu. Petrova \cite{Pet}\footnotemark[16] & $\forall$ & $\forall$ &  $2$ & $q=2$ & \\
\hline 2008 & Oshime \cite{Osh}\footnotemark[17] & $2$ & $0$ &  $\forall$ & $q=\infty$ & \\
\hline 2009 & Watanabe et al \cite{WKNY1} & $3$ & $0$ &  $\forall$ & $q=\infty$ & \\
\hline 2010 & Kalyabin \cite{Kal}\footnotemark[18] & $\forall$ & $0,2$ &  $2$ &  $q=\infty$ & \\
\hline 2010 & \cite{Kal} & $\forall$ & $1$ & $2$ & & $q=\infty$ \\
\hline 2014 & Mukoseeva, Nazarov \cite{MukNaz} & $\forall$ & $4,6$ &  $2$ & $q=\infty$ & \\
\hline 2014 & \cite{MukNaz} & $\forall$ & odd & $2$ & & $q=\infty$ \\
\hline 2021 & Garmanova, Sheipak \cite{GarSh} & $\forall$ & even & $2$ & $q=\infty$ &  \\
\hline 2024 & Garmanova, Sheipak \cite{GarSh_24}\footnotemark[19] & $\forall$ & $n-1$ & $\infty$ & $q=\infty$ &  \\
\hline 2024 & Kazimirov, Sheipak \cite{KazSh} & $\forall$ & $0$ & $\infty$ & $q=\infty$ &  \\
\hline 2024 & Hindov et al \cite{HNOR} & $\forall$ & $0$ & $2$ & $q=1$ &  \\
\hline\end{tabular}
    \caption{The history of inequality (\ref{eq:high})}
    \label{tab3}
\end{table}
\footnotetext{See also \cite{NazPet}.}
\addtocounter{footnote}{1}
\footnotetext{See also \cite{Sl} for the case $n=k+2$, and \cite{MinNaz}.}
\addtocounter{footnote}{1}
\footnotetext{See also \cite{WKNY1}.}
\addtocounter{footnote}{1}
\footnotetext{See also \cite{WKNY2} for the case $k=0$.}
\addtocounter{footnote}{1}
\footnotetext{See also Appendix.}

\begin{conj}
For any $n\ge2$ and $1\le p\le\infty$, the following equality holds:
\begin{equation}
\label{eq:conj2}
\lambda_3(n,1,p,1)=2\lambda_3(n,0,p,\infty).
\end{equation}
Corresponding extremal functions coincide (and are even w.r.t. $x=\frac 12$).
\end{conj}

\begin{remark}
For $n=2$, this statement holds true. Indeed, formulae \eqref{eq:Schmidt}, 
\eqref{eq:Poincare1} and Theorem~\ref{T2} provide $\lambda_3(2,1,p,1)=\frac 14 (p'+1)^{-\frac 1{p'}}$, and \eqref{eq:conj2} follows
from Theorem~\ref{T4}.
\end{remark}

\section{A non-homogeneous inequality for periodic functions}

Finally, we consider an estimate related to the inequality (\ref{eq:Steklov3}):
\begin{equation}
\label{eq:embed}
\|u\|^p_{L^q(0,\ell)}\le \mu\|u\|^p_{W^1_p(0,\ell)}\equiv \mu\int\limits_0^{\ell} \!\big(|u'(x)|^p+|u(x)|^p\big)dx, \quad u(0)=u(\ell).
\end{equation}
By dilation, we can reduce (\ref{eq:embed}) to the case $\ell=1$: 
\begin{equation}
\label{eq:embed1}
\|u\|^p_{L^q(0,1)}\le \widetilde\mu \int\limits_0^1 \big(|u'(x)|^p+m|u(x)|^p\big)dx, \quad u(0)=u(1).
\end{equation}
Here $m>0$, and $\widetilde\mu$ depends on $p$, $q$, and $m$.
\medskip

For $q\le p$, using the H\"older inequality we conclude that the extremal function in (\ref{eq:embed1}) is constant and thus $\widetilde\mu(p,q,m)\equiv m^{-1}$. 

For $q>p$, the problem of sharp constant in (\ref{eq:embed1}) is more delicate and depends on $p$. In \cite{Naz_2000} it was shown that for $p>2$ and any $q>p$ the extremal function in (\ref{eq:embed1}) is non-constant and therefore $\widetilde\mu(p,q,m)> m^{-1}$. In contrast, if $p<2$ then the constant function is a {\bf local} extremal in (\ref{eq:embed1}) for arbitrary $q$. However, given $q>p$, for sufficiently small $m$ the extremal function is constant whereas for sufficiently large $m$ the extremal function is non-constant.\footnote{This problem was considered earlier in \cite{Ka}, but the conclusion in this paper is not correct.}\medskip

The most interesting case is $p=2$. In this case we consider a more general inequality\footnote{For convenience we choose $\ell=2\pi$.} with ``magnetic term''
\begin{equation}
\label{eq:magnetic}
\|u\|^2_{L^q(0,2\pi)}\le \lambda_4 \int\limits_0^{2\pi} \big(|(e^{i\alpha x}u(x))'|^2+m|u(x)|^2\big)dx, \quad u(0)=u(2\pi)
\end{equation}
(here $\lambda_4=\lambda_4(\alpha,m,q)$).

It is easy to see that both sides of (\ref{eq:magnetic}) are invariant w.r.t. replacement $\alpha\mapsto \alpha+k$ and $u(x)\mapsto u(x)e^{-ikx}$, $k\in\mathbb Z$. Therefore, without loss of generality, we may assume that $|\alpha|\le \frac 12$. Then the necessary and sufficient condition of the validity of inequality (\ref{eq:magnetic}) is $m+\alpha^2>0$.
\medskip

Under these assumptions, the following statement holds:

\begin{theorem} 
\label{T10}
Let $\alpha^2(q+2)+m(q-2)\le1$. Then 
$$
\lambda_4(\alpha,m,q)=\frac {(2\pi)^{\frac 2q-1}} {m+\alpha^2}.
$$
The corresponding extremal function is constant.
\medskip

In contrast, if $\alpha^2(q+2)+m(q-2)>1$ then $\lambda_4(\alpha,m,q)>\frac{(2\pi)^{\frac 2q-1}}{m+\alpha^2}$, and the extremal function is non-constant. 
\end{theorem}

\begin{remark}
In the borderline case $\alpha^2(q+2)+m(q-2)=1$ we conclude from $|\alpha|\le\frac 12$
$$
q=\frac {1+2(m-\alpha^2)}{m+\alpha^2}\ge2.
$$
In particular, this means that for $q\le2$ the extremal function in (\ref{eq:magnetic}) is always constant. 
\end{remark}

The history of Theorem \ref{T10} is given in Table \ref{tab4}.

\begin{table}[ht]    
\centering
\begin{tabular} {|c|c|l|}
\hline Year & Authors & Parameters \\
\hline 1999 & Nazarov \cite{Naz_2000}\footnotemark & $\alpha=0$   \\
\hline 2004 & Nazarov \cite{Naz_2013}\footnotemark[23] & $\alpha=0$ \\
\hline 2018 & Nazarov, Shcheglova \cite{NazSch} & $m=0$   \\
\hline 2018 & Dolbeaut et al \cite{DELL} & $\forall \alpha, m$ \\
\hline\end{tabular}
    \caption{The history of inequality (\ref{eq:magnetic})}
    \label{tab4}
\end{table}
\footnotetext{for some values of $m$.}
\addtocounter{footnote}{1}
\footnotetext{In \cite{Naz_2013}, inequalities of arbitrary order were considered. See also \cite{VlKa} for a related result.}

Up to our knowledge, a unique explicit expression for $\lambda_4$ in the non-constancy region of parameters is as follows.

\begin{theorem}[\cite{GaOl}, \S 4, and \cite{GaOl1}]\!\!\footnote{Galunov--Oleinik considered also some higher-order inequalities. See also \cite{ILLZ} for a related result.} \ Let $m+\alpha^2>0$. Then
$$
\lambda_4(m,\alpha,\infty)=
\begin{cases}
\frac 1{2\sqrt{m}}\,\frac {\sinh(2\pi\sqrt{m})}{\cosh(2\pi\sqrt{m})-\cos(2\pi\alpha)}, & m>0;\\
\frac {\pi}{1-\cos(2\pi\alpha)}, & m=0;\\
\frac 1{2\sqrt{-m}}\,\frac {\sin(2\pi\sqrt{-m}\vphantom{\displaystyle 1^1})}{\cos(2\pi\sqrt{-m})-\cos(2\pi\alpha)}, & m<0.
\end{cases}
$$
\end{theorem}

\section{Appendix}

Here we prove that the extremal function in the inequality (\ref{eq:high}) for $p=q=\infty$ and $k=n-1$ is even w.r.t. $x=\frac 12$.

Theorem 2 in \cite{GarSh_24} implies that 
\begin{equation}
\label{L1-approx}
A_{n,n-1,\infty}(a)=\min_P\|\chi_{[0,a]}-P\|_{L^1(0,1)},
\end{equation}
where the mimimum is taken over the set of polynomials of degree no greater than $n-1$.

Let $P_*$ be a minimizing polynomial in (\ref{L1-approx}). Then the necessary condition of minimum reads as follows (see Statement 2 in \cite{GarSh_24} or \cite[Th.~2.1]{Pin}):
\begin{equation}
\label{L1-necess}
\int\limits_0^1 x^{k-1}\,{\rm sign} \big(\chi_{[0,a]}-P_*\big)\,dx=0,\qquad k=1,\dots, n,
\end{equation}
and thus the difference $\chi_{[0,a]}-P_*$ should change sign at $N\ge n$ points (denote these points ${\rm a}_1<{\rm a}_2<\dots<{\rm a}_N$). On the other hand, since the values of $P_*$ at all points ${\rm a}_j$, except maybe for ${\rm a}_j=a$, are prescribed, we have $N\le n+1$.

Let $P$ be another minimizing polynomial.\footnote{Notice that $L_1(0,1)$ is not a strictly convex space, so in general we have no uniquess of the best approximation function.} Then Proposition 2.4 in \cite{Pin} shows that
$$
(\chi_{[0,a]}-P_*)(\chi_{[0,a]}-P)\ge0 \qquad \mbox{a.e. on } [0,1].
$$
This means that the difference $\chi_{[0,a]}-P$ changes sign at the same points. 

Statement 3 in \cite{GarSh_24} implies that the function $w_{n,a}$ providing the value $A_{n,n-1,\infty}(a)$ is the $n$-th primitive of the function ${\rm sign} \big(\chi_{[0,a]}-P_*\big)$. So, if we choose $a$ as the maximum point of $A_{n,n-1,\infty}$ given in Theorem \ref{T7}, then $w_{n,a}^{(n-1)}$ reaches its maximum at $a$, and therefore $a$ belongs to the set $\{{\rm a}_j\}_{j=1}^N$.

The condition (\ref{L1-necess}) can be rewritten as follows:
\begin{equation}
\label{a_j-eq}
{\rm a}_1^k-{\rm a}_2^k+
\ldots
+(-1)^{N-1}{\rm a}_N^k+
\frac {(-1)^N}2=0,
\qquad k=1,\ldots, n.
\end{equation}
It is easy to prove (cf. \cite[Cor.~13]{CHJSV}) that
\begin{equation}
\label{a_j-sol}
{\rm a}_j=\mathfrak{a}_j=\sin^2\tfrac{\pi j}{4(n+1)}, \qquad j=1,\dots,n,
\end{equation}
satisfy the equations (\ref{a_j-eq}) with $N=n$. Notice that the maximum point of $A_{n,n-1,\infty}$ belongs to the set $\{{\rm a}_j\}_{j=1}^n$, as required.
 
It remains to observe that equalities (\ref{a_j-sol}) immediately imply symmetry of the function ${\rm sign} \big(\chi_{[0,a]}-P_*\big)$: it is even or odd w.r.t. $x=\frac 12$ if $n$ is even or odd, respectively. Since the extremal function in the original inequality is its $n$-th primitive, it is even w.r.t. $x=\frac 12$ in both cases, and we are done.\bigskip

{\bf Acknowledgements.} We are grateful to Prof. N.G. Kuznetsov for valuable comments and suggestions, and to Prof. I.A. Sheipak who provided us with some important references.

The work of the first author was supported by the Ministry of Science and Higher Education of the Russian Federation (agreement 075-15-2025-344 dated 29/04/2025 for Saint Petersburg Leonhard Euler International Mathematical Institute).

%\small

\end{document}